\documentclass[12pt]{amsart}
\usepackage[francais]{babel}
\usepackage{graphics}
\usepackage{graphicx}
\usepackage{pdfpages}

\title{Miroir circulaire et polyn\^omes de Stewart} 

\author{Jean-Claude Carrega et Labib Haddad}
\address{Jean-Claude Carr\Žga, R\Žsidence Horizon,
12 Bd de l'Europe, 69110 Sainte-Foy-l\s-Lyon, FRANCE.}
\email{jeanclaudecarrega@orange.fr}
\address{Labib Haddad, 120 rue de Charonne,75011 Paris, FRANCE.}
\email{labib.haddad@wanadoo.fr}

\usepackage{amssymb}
\usepackage{amsmath}
\usepackage{endnotes}

\newcommand{\su}{\subsection*}
\newcommand{\head}{\section*}
\newcommand{\noi}{\noindent}

\newcommand{\Ž}{\'e}

\newcommand{\ˆ}{\`a}
\newcommand{\}{\`u}

\newcommand{\A}{\mathbb A}

\newcommand{\Q}{\mathbb Q}
\newcommand{\R}{\mathbb R}

\newcommand{\geqs}{\geqslant}

\newcommand{\ali} {\begin{aligned}}   
\newcommand{\ala} {\end{aligned}}

\newcommand {\et}{\ \text{et}\  }

\newcommand {\ou}{\ \text{ou}\  }

\newcommand{\bc}{\begin{cases}}
\newcommand{\ec}{\end{cases}}

\begin{document}
\maketitle

\thispagestyle{empty}

\

Comment peuvent se  conjuguer deux variations sur un m\me th\me, celui des constructions g\Žom\Žtriques \ˆ l'aide de la r\gle et du compas,  tel est le sujet de cette petite note.

\su{Le polyn\™me de Stewart} {\sl Tout nombre constructible  est alg\Žbrique sur le corps $\Q$ et son degr\Ž est une puissance de $2$}. On le sait. On sait \Žgalement que la r\Žciproque est fausse : il y a des nombres non constructibles qui sont alg\Žbriques sur $\Q$ et dont le degr\Ž est une puissance de $2$. Le premier auteur de cette note en a  donn\Ž les d\Žmonstrations dans un ouvrage destin\Ž aux \Žtudiants et aux enseignants du  secondaire, {\sl Th\Žorie des corps, La r\gle et le compas} [1]. 
 
\ 
 
 \noi Il a pr\Žsent\Ž, en particulier, le polyn\™me $$X^4 - X -1$$
comme contre-exemple \ˆ la r\Žciproque. Il \Žtablit que ce polyn\™me est irr\Žductible sur $\Q$ et que l'une de ses deux racines r\Želles n'est pas constructibles \ˆ l'aide de la r\gle et du compas bien qu'elle
soit alg\Žbrique de degr\Ž $4 = 2^2$. Il attribue ce contre-exemple \ˆ {\sc Ian STEWART}. En fait, le groupe de Galois de ce polyn\™me est le groupe sym\Žtrique $S_4$ (comme on le voit, d'un clic, \ˆ l'aide du logiciel Maple) et, bien entendu, aucune de ses deux racines r\Želles n'est constructible.

\

\noi En effet, dans le langage d'aujourd'hui, on \Žnonce comme suit la caract\Žrisation des nombres constructibles :   un nombre $\alpha$ est constructible, si et seulement si l'ordre de son  groupe de Galois associ\Ž, $\mathbf G(\alpha)$, est une puissance de $2$, ce qui n'est pas le cas du groupe sym\Žtrique $S_4$. Le r\Žsultat est d\Žtaill\Ž dans le livre [1] d\Žj\ˆ cit\Ž. On le trouve bien entendu, \Žgalement, dans le van der Waerden [4]. On sait aussi que l'ordre du groupe de Galois $\mathbf G(\alpha)$ est le degr\Ž Ê$[K : \Q]$  de l'extension  $K$ de $\Q$, o\ $K$ est le  corps de d\Žcomposition du polyn\™me minimal de $\alpha$. Voir [1], par exemple.

\includegraphics[width=6.5in]{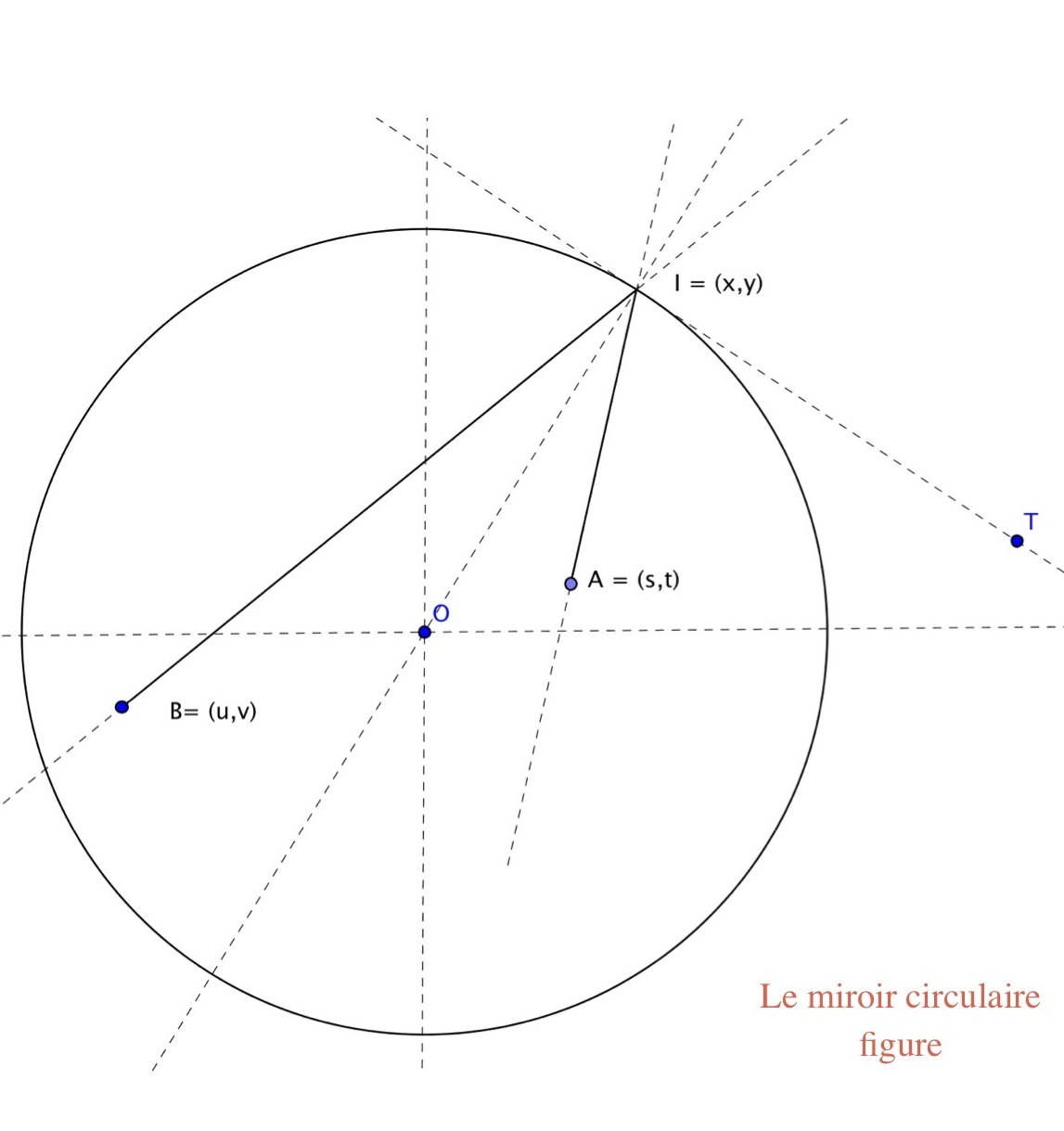}
\includegraphics[height=21in]{f4}
\includegraphics[scale=3]{f4}

\su{Le miroir circulaire} Dans le plan, on se donne  une circonf\Žrence, $C$, et deux points, $A$ et $B$. On cherche \ˆ construire les points $I$ de $C$ en lesquels le rayon lumineux $AI$ se r\Žfl\Žchit pour repasser  par $B$. Voir la figure ci-dessus. 

\

\noi C'est ce que l'on appelle le probl\me du miroir circulaire. On sait que cette construction est impossible \ˆ l'aide de la r\gle et du compas, sauf dans des cas tr\s particuliers, pour certaines positions {\it critiques} des points $A$ et $B$.

\

\noi On a, parfois, appel\Ž {\it probl\me d'Alhazen} ce probl\me du miroir circulaire (ou sph\Žrique), voir  [3]. On l'appelle aussi {\it probl\me du billard circulaire} (voir {\sc Carrega} [1,page 253, exercice 25]).

\

Voici le r\Žcit d'un lien que l'on tisse entre le miroir circulaire et  le polyn\™me de Stewart. 

\

\head{Les polyn\™mes de Stewart}

\

Plus g\Žn\Žralement, appellons {\it polyn\™me de Stewart} tout polyn\™me de la forme
$$S(X) = X^4 - rX -1 \ \ \text{o\ $r$ est un nombre rationnel non nul},$$
ainsi que tous leurs multiples scalaires  $\lambda S(X)$ o\ $\lambda$ est un nombre rationnel non  nul.

\

\noi Quitte \ˆ changer $X$ en $-X$, on peut se ramener aux cas o\ l'on a  $r > 0$. 

\

\noi En imitant ce qui est fait dans {\sc Carrega} [1, page 39], partant du polyn\™me
$$S(X) = X^4 - rX-1,$$
on pose 
$$r = a\sqrt{a^4+4} \ , \ 2b = a^2 + \sqrt{a^4+4} \ , \ 2\bar b = a^2 - \sqrt{a^4+4}.$$

Il vient :
$$S(X) =  (X^2 + aX + b)(X^2-aX + \bar b) \ , \ a^6 + 4a^2-r^2 = 0.$$
Ainsi, $a^2$ est racine du polyn\™me  $R(Y) = Y^3 + 4Y - r^2$ de degr\Ž $3$. Ce polyn\™me  est fonction strictement croissante de  $Y$; il prend des valeurs n\Žgatives pour $Y < 0$ et poss\de donc une seule racine r\Želle laquelle est $>0$. \`A cette racine r\Želle, $a^2$, correspondent deux valeurs oppos\Žes $\pm a$. On doit prendre, bien \Žvidemment, la valeur qui a le m\me signe que $r$. On notera ainsi que $a,b \et \bar b$ sont des nombres r\Žels, alg\Žbriques. Le discriminant du polyn\™me $X^2 + aX + b$ est \Žgal \ˆ $a^2 - 4b$, i.e., $-2\sqrt{a^4+4} -a^2<0$; 
celui du polyn\™me {\it conjugu\Ž} $X^2-aX + \bar b$ est \Žgal \ˆ $2\sqrt{a^4+4} -a^2>0$. 

\

\noi Le polyn\™me de Stewart $S(X)$ a toujours, ainsi, $4$ racines distinctes, soit deux racines complexes conjugu\Žes et deux racines r\Želles, \ˆ savoir
$$\frac{-a \pm i\sqrt{2\sqrt{a^4+4} + a^2}}{2} \ , \ \frac{a \pm \sqrt{2\sqrt{a^4+4} - a^2}}{2}.$$

\

\su{Irr\Žductibilit\Ž} Soit $\A$ le corps  des nombres alg\Žbriques r\Žels. Le polyn\™me de Stewart $S(X)$ est ainsi d\Žcomposable dans $\A[X]$ en un produit de deux polyn\™mes de degr\Ž $2$  : le polyn\™me $U(X) = X^2 + aX + b$ et le polyn\™me $V(X) = X^2-aX + \bar b$. Pour qu'il soit d\Žcomposable  dans $\Q[X]$, il faut et il suffit que l'une des deux conditions suivantes soit remplie : 

\

\noi (C1) Il existe un nombre rationnel $s$ non nul tel que l'on ait $r = s^3-1/s$.

\noi (C2) Il existe un nombre rationnel $a$ tel que l'on ait $r^2 = a^6 + 4a^2$.

\

\noi {\bf En effet}, la condition C1 revient \ˆ dire que le polyn\™me $S(X)$ poss\de une racine rationnelle ce qui \Žquivaut \ˆ dire qu'il se d\Žcompose dans $\Q[X]$ en un produit d'un polyn\™me de degr\Ž $1$ avec un polyn\™me de degr\Ž $3$. Quant \ˆ la condition C2, elle implique que $\sqrt{a^4 + 4} = r/a$ est un nombre rationnel. Cel\ˆ entra\"ne que les coefficients $a, b \et \bar b$ sont rationnels, de sorte que les polyn\™mes $U(X)$ et $V(X)$ appartiennent \ˆ $\Q[X]$ : le polyn\™me $S(X)$ se d\Žcomposerait alors dans $\Q[X]$ en un produit de deux polyn\™mes du second degr\Ž. R\Žciproquement, si $S(X)$ \Žtait produit de deux polyn\™mes de $\Q[X]$, du second degr\Ž, ces deux polyn\™mes seraient $U(X)$ et $V(X)$, {\it par n\Žcessit\Ž}, d'o\ le r\Žsultat. \hfill {\bf cqfd}

\

\su{Nota} On signale toutefois ceci. La condition C2 n'est jamais satisfaite : on en donnera la d\Žmonstration dans l'{\sc Appendice} ci-dessous. Cela veut dire q'un polyn\™me de Stewart ne se d\Žcompose jamais en un produit de deux polyn\™mes du second degr\Ž dans $\Q[X]$.  Autrement dit, $S(X)$ est r\Žductible sur $\Q$ si et seulement s'il poss\de une racine rationnelle. 
 
\

\noi Dans tous les autres cas, le polyn\™me $S(X)$ est irr\Žductible sur $\Q$, en particulier pour $r=1$, le cas {\it princeps}.

\

\su{Le cas o\ $\mathbf r$ est entier} Plus g\Žn\Žralement, lorsque $r$ est un entier (non nul), le polyn\™me  $S(X)$ est irr\Žductible. Pour le voir, on montre que la condition C1 n'est pas satisfaite. Si l'on avait  $r = s^3 - 1/s$ o\  $s=p/q$ est une fraction irr\Žductible, on aurait $p^4  - rpq^3 - q^4 = 0$, donc $p$ diviserait $q$ et $q$ diviserait $p$,  de sorte que $s=p/q$ serait \Žgale \ˆ $\pm 1$ et $r$ serait nul, ce qui n'est pas ! 

\

\noi Lorsque $r$ est entier, si $a^2$ est rationnel, c'est un entier. En effet, $a^2$ est  racine r\Želle positive du polyn\™me $R(Y)$. Si $a^2$ est rationnelle,  on l'\Žcrit sous forme irr\Žductible $a^2 = p/q$. Il vient
$p^3	+ 4 pq^2	- r^2	q^3	= 0$, donc $q$ divise $p$ et $a^2$ est un entier. 

\

\noi Lorsque $r$ est entier, on peut montrer simplement, directement, que la condition C2 n'est pas satisfaite, comme suit. Si l'on avait  $r = a\sqrt{a^4+4}$ pour $a$ rationnel, $a^2$ serait entier, d'apr\s ce qui pr\Žc\de. On aurait $a^4+ 4 = (r/a)^2$, de sorte que $u= r/a$ serait entier et l'on aurait
$$4 = u^2 - a^4 = (u-a^2)(u+a^2) \ , \ 0 < u-a^2 < u+a^2.$$
La seule possibilit\Ž serait alors $u-a^2 = 1$ et $u+a^2 =4$ ce qui entra\"ne
$2u = 5$, impossible puisque $u$ est entier.

\

\su{Le cas o\ $\mathbf r$ est un nombre premier} On suppose que $r$ est un nombre {\bf premier}. On  \Žtablit que le polyn\™me $R(Y)$ est irr\Žductible. Pour cela,  on doit montrer que sa racine $t=a^2$ n'est pas rationnelle. D'apr\s ce qui pr\Žc\de, il suffit de montrer que $t$ n'est pas un entier !

\su{D\Žmonstration} On a $t(t^2+4) = r^2$. Si $t$ \Žtait entier, il diviserait $r^2$. Or, $r$ \Žtant premier, on ne peut avoir que $t=1 \ou t=r \ou t = r^2$.

\noi Si $t=1$, on aurait $r^2 = 5$ qui est impossible.

\noi Si $t = r$, on aurait  $r^2 -r +4 = 0$ qui n'a pas de racine r\Želle.

\noi Si $t = r^2$, on aurait $r^4 + 4 = 1$ qui est impossible.\qed

\

\noi Ainsi $t= a^2$ n'est pas un nombre constructible et il en r\Žsulte que $a$ non plus n'est pas constructible. Les deux racines r\Želles du second facteur $X^2 - a X + \bar b$ de $S(X)$  ont pour somme $a$, donc l'une au moins de ses racines n'est pas constructible. Ainsi, l'ordre du groupe de Galois de $S(X)$ n'est pas une puissance de deux.

\

\noi En fait, aucune des $2$ racines r\Želles de ce polyn\™me n'est constructible car toutes deux ont le m\me polyn\™me minimal $S(X)$ [et l'ordre du groupe de Galois de $S(X)$ n'est pas une puissance de deux].

\

\noi  Ainsi, les polyn\™mes de Stewart $S(X)$ avec $r$ nombre premier permettent d'obtenir, par leurs racines r\Želles, une infinit\Ž de nombres alg\Žbriques de degr\Ž $4$ qui ne sont pas constructibles.Ê

\

\su{Une variante}  \`A la fin de la d\Žmonstration pr\Žc\Ždente, on peut tout aussi bien utiliser l'exercice 24 de [1],  pages 252-253. Dans le a) de cet exercice, on fait d\Žmontrer le r\Žsultat g\Žn\Žral suivant au sujet des polyn\™mes  de degr\Ž $4$ : \sl une racine r\Želle d'un polyn\™me irr\Žductible $P(X)\in \Q[X]$ de degr\Ž $4$ est constructible si et seulement si le r\Žsolvant de $P(X)$ est r\Žductible sur $\Q$\rm. En l'occurence, le r\Žsolvant du polyn\™me de Stewart $S(X)$ n'est autre que le polyn\™me $Y^3 + 4 Y + r^2= -R(-Y)$, ce qui ach\ve la d\Žmonstration.

\su{Remarque} Examinons l'exemple o\ $r=4$, un entier non premier. Dans ce cas, $2$ est racine de $R(Y) = Y^3 + 4Y - 16$, de sorte que l'on a $t = a^2 = 2$ et $a = \sqrt 2$. Les $2$ racines r\Želles de $S(X)$, donn\Žes par les formules ci-dessus (en haut de la page 4)  avec $a = \sqrt 2$, sont constructibles et alg\Žbriques et de degr\Ž $4$.

\

\head{Le miroir circulaire}

\

Dans le plan, on se donne  une circonf\Žrence, $C$, et deux points, $A$ et $B$. On cherche les points $I$ de $C$ en lesquels le rayon lumineux $AI$ se r\Žfl\Žchit pour repasser  par $B$.

\

Voici une solution analytique, \ˆ suivre sur la figure en page 2,
ci-dessus. 

\

\noi Dans le plan des $x,y$, ayant O pour origine, on prend la circonf\Žrence
\[x^2 + y^2 = 1\tag{C}\]
ainsi que les points $A=(s,t)$ et $B = (u,v)$. Soit $IT$ la tangente en $I$ \ˆ $C$. On voudrait trouver les points $I=(x,y)$ de $C$ tels que les droites $IO$ et $IT$ soient les deux bissectrices des angles que forment les droites 
$IA$ et $IB$. 

\

\noi Il faut et il suffit pour cela que le rapport anharmonique des pentes des droites $IA, IB, IO, IT,$ soit \Žgal \ˆ $-1$. On a
$$p_{IA} = \frac{y-t}{x-s} \ , \ p_{IB} = \frac{y-v}{x-u} \ , \ p_{IO} = \frac{y}{x} \ , \ p_{IT} = -\frac{x}{y}.$$
Ainsi, il faut et il suffit que l'on ait :
$$\frac{\dfrac{y-t}{x-s}- \dfrac{y}{x}}{\dfrac{y-v}{x-u}-  \dfrac{y}{x}}: \frac{\dfrac{y-t}{x-s}+ \dfrac{x}{y}}{\dfrac{y-v}{x-u}+  \dfrac{x}{y}}= -1.$$
Tous calculs faits, cela donne :
$$(tx - sy)(x^2+y^2-ux -vy) + (vx-uy)(x^2+y^2 -sx -ty) =0.$$
Puisque $x^2 + y^2 = 1$, on obtient
\[(sv+tu)(y^2-x^2) + 2(su -tv)xy   + (t+v)x - (s+u)y= 0,\tag{H}\]
l'\Žquation d'une hyperbole, H, passant par l'origine O.

\

\noi En utilisant la param\Žtrisation classique suivante du cercle  
 $$x = \frac{1-z^2}{1+z^2} \ , \ y = \frac{2z}{1+z^2},$$
et tous calculs faits,  l'\Žquation $H$ prend la forme que voici :
$$\frac{Q(z)}{(1+z^2)^2} = 0,$$
o\ $Q(z)$ est le polyn\™me suivant de degr\Ž $4$ en $z$ :

\

\fbox{\begin{minipage}[t]{12cm}
$$Q(z)= (sv+tu + t+v)z^4 + 2( 2su- 2tv+s+ u)z^3$$
$$- 6(sv+tu)z^2 - 2(2su- 2tv-s-u)z + (sv+tu -t-v).$$

\end{minipage}}

\

\

\noi Par commodit\Ž, on dira que ces polyn\™mes $Q(z)$ ainsi que tous leurs multiples scalaires $\lambda Q(z)$ sont  les {\it polyn\™mes d'Alhazen}.

\

\noi On montre alors que tout polyn\™me de Stewart est un polyn\™me d'Alhazen.

\

\noi Pour cela, on sp\Žcialise une premi\re fois, en prenant
$$t=s \et v=-u.$$
Le polyn\™me $Q(z)$ prend la forme
$$Q(z)= (s-u)z^4+ 2(4su+s+u)z^3- 2(4su-s-u)z-(s-u).$$
On sp\Žcialise de nouveau, en prenant
$$u = \frac{-s}{4s+1}.$$
Le polyn\™me $Q(z)$ s'\Žcrit :
$$Q(z) = \frac{2s(2s+1)}{4s+1}(z^4-rz-1) \ \text{o\} \ r= \frac{-8s}{2s+1}.$$
Il en d\Žcoule, comme annonc\Ž, que {\bf tout polyn\™me de Stewart est un polyn\™me d'Alhazen}.

\

\noi Le polyn\™me de Stewart $S(X) = X^4-rX-1$ d\Žpend du seul param\tre $r$. Le polyn\™me d'Alhazen $Q(z)$ d\Žpend des $4$ param\tres $s,t,u,v,$ les coordonn\Žes des points $A$ et $B$. 

\

\noi Dans {\sc Carrega} [1, p. 266, solution de l'exercice 25 sur le billard circulaire]  on trouve l'expression suivante du polyn\™me $Q(z)$
$$(a+1)cz^4 +2(a+b+2ab)z^3 -6acz^2 + 2(a+b-2ab)z + (a-1)c,$$
obtenue en utilisant les nombres complexes pour exprimer l'\Žgalit\Ž des deux arguments correspondant aux angles d\Žfinis par la bissectrice,
cela \Žtant fait dans le cas particulier o\ $s = a, t=0, u=b, v=c$.

\

\centerline{$*$ \ \ $*$ \ \ $*$}

\

Cette \Žtude a permis la rencontre improbable des noms de deux math\Žmaticiens que 10 si\cles s\Žparent :
Ian Stewart est professeur \Žm\Žrite \ˆ l'universit\Ž de Warwick en Angleterre. Il est l'auteur de nombreux ouvrages remarquables.
Ibn Al Haytham (965 - 1039), connu en Occident sous le nom de Alhazen, est un savant du monde m\Ždi\Žval arabo-musulman, originaire de Perse. Il est l'auteur de trait\Žs sur la G\Žom\Žtrie, l'Optique et l'Astronomie.

\

\centerline{$*$ \ \ $*$ \ \ $*$}

\

\noi Pour les groupes de Galois des \Žquations de degr\Ž 3 et 4, on pourra consulter utilement le livre de Kaplansky, [2].

\

\

\

\

\

\head{Appendice}

\

\

\

{\sc Euler} a montr\Ž ceci : {\sl La somme de deux bicarr\Žs d'entiers non nuls n'est jamais le carr\Ž d'un entier non nul}. Autrement dit,  l'\Žquation $x^4 + y^4 = z^2$ n'a pas de solutions en entiers $x, y, z,$ strictement positifs. 

\

\noi \{Voir \ˆ ce sujet le livre de \sc L. E. Dickson\rm, \sl Theory of numbers\rm, vol. II, pages 615 et s. o\ on pourra lire la longue histoire de l'\Žquation $x^4 + y^4 = z^4$. On pourra \Žgalement consulter le livre de {\sc Pierre SAMUEL}, \sl Th\Žorie alg\Žbrique des nombres\rm, Collection M\Žthodes, Hermann, Paris 1967, Deuxi\me \Ždition revue et corrig\Že, Paris, 1971,  page 21. On y trouve aussi, page 20, la r\gle de Diophante dont il sera question ci-dessous.\}

\

\noi Plus g\Žn\Žralement, on a le r\Žsultat suivant, lequel est un cas tr\s particulier du Th\Žor\me 169 de {\sc Hilbert}, bien plus g\Žn\Žral. On pourra consulter \rm  Th\Žorie des corps de nombres algŽbriques, deuxi\me partie, Trad. A. LEVY, \it Annales de la facult\Ž des sciences de Toulouse  $3^e$ s\Žrie\rm, 
tome 2, \no 3-4, p. 455-456.

\su{Th\Žor\me} \sl L'\Žquation diophantienne 
$x^4 + 4y^4 = z^2$ n'a pas de solutions en nombres entiers strictement positifs\rm.

\su{D\Žmonstration} On se servira du r\Žsultat suivant connu sous le nom de \it r\gle \rm  de Diophante. \sl Les solutions de l'\Žquation $x^2 + y^2 = z^2$ en $x, y, z,$ entiers strictement positifs  et premiers entre eux, sont de la forme 
$$x = a^2 - b^2 \ , \ y = 2ab \ , \ z = a^2 + b^2,$$
o\ $a$ et $b$ sont entiers, strictement positifs, premiers entre eux, l'un pair et l'autre impair\rm. 

\

\noi On utilise la m\Žthode de la descente,  en supposant que l'\Žquation diophantienne $x^4 + 4y^4 = z^2$ poss\de des solutions en nombres entiers strictement positifs. On se donne une des solution, $(x,y,z)$,  pour laquelle $z$ est le plus petit possible. On observe que $x, 2y, z,$ sont alors deux \ˆ deux premiers entre eux. En effet, si un nombre premier impair $p$ divise deux d'entre eux, il divise le troisi\me et l'on aurait 

\

$x = pu \ , \ y=pv \ , \ z= pg,$

$p^4u^4 + 4p^4v^4 = p^2g^2,$

$p^2(u^4 + 4v^4) = g^2,$

$g$ serait divisble par $p$ et l'on aurait $g = pw$, d'o\

$u^4 + 4 v^4 = w^2$ o\ $w < z$, ce qui est impossible.

\

\noi De m\me, si $2$ divisait $x$ ou $z$, il diviserait  les deux et l'on aurait

\

$x = 2u \ , z = 2w$,

$16u^4 + 4y^4 = 4w^2$,

$4u^4 + y^4 = w^2$ o\ $w< z$, ce qui est \Žgalement impossible.

\

\noi On \Žcrit $(x^2)^2 + (2y^2)^2 = z^2$. En vertu de la r\gle de Diophante, on aurait

\ 

$x^2 = a^2 - b^2 \ , \ 2y^2 = 2ab \ , \ z = a^2 + b^2$, 

$x^2 = a^2 - b^2 \ , \ y^2 = ab \ , \ z = a^2 + b^2$, 

\

\noi o\ $a$ et $b$ sont des entiers non nuls, premiers entre eux, l'un pair et l'autre impair. Mais alors $x$ est impair donc $a^2 - b^2= x^2 \equiv 1 \mod 4$,  de sorte que $a$ est impair et $b$ est pair !

\

\noi Or, $ab = y^2$ est un carr\Ž, donc $a$ et $b$ sont des carr\Žs. De plus, on a $x^2 +b^2 = a^2$, donc $a = m^2 + n^2 \et b = 2mn$ o\ $m$ et $n$ sont des entiers non nuls, premiers entre eux, l'un pair et l'autre impair (par la r\gle de Diophante). Sans nuire \ˆ la g\Žn\Žralit\Ž, on peut supposer que c'est $m$ qui est pair.

\

\noi Puisque $b = (2m)n$ est un carr\Ž, $2m$ et $n$ sont des carr\Žs. On aurait ainsi $2m = 4u^2 \et n=v^2$ Ainsi $4u^4 + v^4 = m^2 + n^2 = a$ o\ $a$ est un carr\Ž. Or,  $a = \sqrt{z- b^2} < z^2$, autrement dit $\sqrt a < z$,  ce qui est impossible.\qed

\

\noi De cela, on d\Žduit ais\Žment ceci : pour $r$ rationnel non nul, il n'exite pas de nombre rationnel $a$ tel que l'on ait $(r/a)^2 = a^4 + 4$ car, en \Žcrivant $a = y/z$ comme fraction irr\Žductible, on aurait 
$(z^2r/a)^2=y^4 +4z^4$ o\ $x = z^2r/a$ serait entier, ce qui est impossible.

\

\noi Autrement dit, comme annonc\Ž, la condition C2 n'est jamais satisfaite.

\

\

\

\

\head{Bibliographie}

\

1. Jean-Claude CARREGA, {\it Th\Žorie des corps, La r\gle et le compas},  Nouvelle \Ždition enrichie d'exercices,  Collection Formation des enseignants, Hermann, Paris, 1989.

\

2.  Irvin KAPLANSKY,   {\it Fields and rings}, U. Chicago Press, (en particulier, p.50-52).

\

3. Peter M. NEUMANN,  {\it Reflections on reflection in a spherical mirror},
Amer. Math. Monthly, {\bf105} (1998) No. 6, 523-528.

\

4. van der WAERDEN, {\it Modern Algebra}, Tome 1, p.183-187,  \S 59, [dans l'\Ždition  Frederik Ungar Publishing Co., 1949.]

\

\

\

\

\enddocument